\documentclass[12pt]{amsart}
\usepackage{amssymb}
\usepackage{amscd}

\DeclareMathOperator{\Ad}{Ad}
\DeclareMathOperator{\ad}{ad}

\DeclareMathOperator{\im}{Im}

\newcommand{\G}{\mathcal{G}}
\renewcommand{\L}{\mathcal{L}}
\newcommand\lie[1]{\mathfrak{#1}}
\renewcommand{\d}{\lie{d}}
\newcommand{\g}{\lie{g}}
\newcommand{\h}{\lie{h}}

\newcommand{\C}{\mathbb{C}}
\newcommand{\R}{\mathbb{R}}

\renewcommand{\l}{\langle}
\renewcommand{\r}{\rangle}

\newcommand{\e}{\varepsilon}

\numberwithin{equation}{subsection}
\numberwithin{enumi}{equation}

\theoremstyle{plain}
\newtheorem{theorem}[equation]{Theorem}
\newtheorem{proposition}[equation]{Proposition}

\theoremstyle{definition}
\newtheorem{definition}[equation]{Definition}

\theoremstyle{remark}
\newtheorem{remark}[equation]{Remark}

\theoremstyle{definition}
\newtheorem{example}[equation]{Example}

\begin{document}

\title{Manin pairs and moment maps}

\author{Anton Alekseev}
\address{Institute for Theoretical Physics 
\\ Uppsala University \\
Box 803 \\ \mbox{S-75108} Uppsala \\ Sweden}
\email{alekseev@teorfys.uu.se}

\author{Yvette Kosmann-Schwarzbach}
\address{U.M.R. 7640 du C.N.R.S. \\ Centre de Math{\'e}matiques \\
Ecole Polytechnique \\ \mbox{F-91128} Palaiseau \\ France}
\email{yks@math.polytechnique.fr}

\date{}

\begin{abstract}
A Lie group $G$ in a group pair $(D,G)$, integrating a Lie algebra 
${\mathfrak g}$
in a Manin pair $({\mathfrak d}, {\mathfrak g})$, 
has a quasi-Poisson structure. We define the quasi-Poisson actions of
such Lie groups $G$, that generalize the Poisson actions of Poisson Lie
groups.
We define and study the moment maps for those quasi-Poisson actions
which are quasi-hamiltonian. These moment maps take
values in the homogeneous space $D/G$.
We prove an analogue of the hamiltonian reduction theorem for
quasi-Poisson group actions, and we study the 
symplectic leaves of the orbit spaces of quasi-hamiltonian spaces.
\end{abstract}


\maketitle


\section{Introduction}

The purpose of this article is to provide a framework for
Lie-group valued moment map theories. In the
usual theory (see, {\it e.g.}, \cite{GS}), 
the moment map corresponding to an action of a Lie group $G$ 
on a symplectic manifold $(M,\omega)$ takes values in the
dual space $\g^*$ of the Lie algebra $\g$. In the case
of $G=S^1$, an abelian Lie-group valued moment map taking values 
in $S^{1}$ instead of $u(1)^{*}=\R$ was introduced in \cite{McD}  
by McDuff.
The first nonabelian Lie-group valued moment map theory
to be proposed
was that of Lu and Weinstein \cite{Lu}, \cite{LW}.
In their approach, $G$  is a Poisson Lie group,  
$M$ has a symplectic form $\omega$ (or, more
generally, a Poisson bivector)
which is not $G$-invariant,
and the moment map takes values in the 
Poisson Lie group dual to $G$.
When $G$ is a compact semi-simple Lie group, 
the target space is
the symmetric space 
$G^{\C}/G$ which is equipped with the structure of a Lie group.
Recently, another nonabelian moment map theory where the moment map
takes values in the same compact simple Lie group $G$ as the
one which acts on the manifold has been developed \cite{AMM1}.
In this theory, the 2-form $\omega$ is $G$-invariant but not closed.

These various formulations of the moment map theories 
share many features.
For instance, they all have a well defined
notion of a hamiltonian reduction \cite{Lu}, \cite{AMM1},
convexity properties \cite{FR}, \cite{A}, \cite{MW} and
localization formulas \cite{G}, \cite{W}, \cite{AMW}. We think that these 
common features justify a search for a unified formulation.

Our proposal of a general moment map theory is based on the notions
of a Manin pair and of a Manin quasi-triple, and our
main technical tool is the theory of quasi-Poisson Lie groups
developed in \cite{YKS1} and \cite{YKS2}. We introduce the notion of
a quasi-Poisson space as a space with a $G$-action and a bivector
$P$ such that the Schouten bracket $[P,P]$ is expressed as
a certain trilinear combination of the vector fields generating
the $G$-action. These spaces are the natural objects 
upon which quasi-Poisson Lie groups act. 
We define and study the moment maps 
for the actions on a quasi-Poisson space 
of a Lie group $G$ whose Lie algebra, $\g$, belongs to a 
Manin pair $(\d,\g)$. 
Our formulation is close
in spirit to the one of Lu and Weinstein, which is based on the 
notions of a Manin triple and of a Poisson Lie group.
While Manin triples and Poisson Lie groups had been introduced as the 
classical 
limits of quantum groups \cite{Dr1}, the generalized objects that we
consider here are classical limits of quasi-Hopf algebras \cite{Dr2}. 
The present generalization to Manin quasi-triples and quasi-Poisson Lie
groups is necessary in order to provide a conceptual
explanation for the theory of group valued moment maps introduced 
and developed in \cite{AMM1} and \cite{AMW}.

Section 2 includes some basic information on Manin pairs and
Manin quasi-triples. In Section 3 we introduce ``group pairs'' 
which integrate Manin pairs. (Group pairs with a 
symmetric structure were also studied in \cite{Le}.)
The definition of
quasi-Poisson actions is presented in Section 4. If $(D,G)$ is a
group pair, then the dressing action of $G$ on $D/G$ is quasi-Poisson.
We show that a given
action remains quasi-Poisson when both the quasi-Poisson Lie group and
the quasi-Poisson space on which it acts are modified by a twist.
Thus the notion of a quasi-Poisson action is related to that of a Manin
pair, rather than to a particular Manin quasi-triple.
We prove that if
$(M,P)$ is a quasi-Poisson space with a bivector $P$ acted upon by
a quasi-Poisson Lie group, the bivector $P$ induces a genuine Poisson
bracket on the space of $G$-invariant smooth functions on $M$. 
The definition of generalized moment maps is given in
Section 5. In our setting, the moment maps are always assumed
to be equivariant, a property which the usual moment maps 
may or may not have. We then prove that moment maps
are bivector maps. This generalizes the theorem stating that the moment
map is equivariant if and only if it is a Poisson map.
The dressing action of $G$ on $D/G$ has the identity of $D/G$
as a moment map.
We show that, on the quasi-Poisson spaces that admit a moment map 
-- called quasi-hamiltonian spaces --, 
there is a well-defined distribution 
defining a generalized foliation by 
nondegenerate quasi-hamiltonian spaces that contain the $G$-orbits.
We also prove that, in a quasi-Poisson space with a nondegenerate
bivector $P$, the symplectic leaves in the orbit space are
connected components of the projection of the
level sets of the moment map. This is a generalization of the usual
symplectic reduction theorem.

The fundamental example of a Manin pair is that, described in
Example \ref{gminus}, of a Lie algebra $\g$ 
with an invariant scalar product,
embedded diagonally in the double $\g \oplus \g$, 
equipped with the corresponding
hyperbolic metric. For the corresponding group pair $(G\times G,G)$,
the target of the moment map $(G \times G) / G$ can be identified with $G$
itself. It can be shown that any hamiltonian $G$-space
with group-valued moment map in
the sense of \cite{AMM1} and \cite{AMW} 
is also a quasi-hamiltonian space associated to this group pair, in the
sense of Definition \ref{momentdef}.
This can be easily shown (see Section 5.2) in the case of abelian Lie
groups and will be proved in the general case in a subsequent paper.

\section{Manin pairs and $r$-matrices}

In this section we introduce the notions of a Manin pair and of a classical
$r$-matrix which provide the main building blocks of the generalized
moment map  theory. All vector spaces are over $\R$ or $\C$ and, for
simplicity, we assume that they are finite-dimensional.
We use the Einstein summation convention.

\subsection{Manin pairs and Manin quasi-triples}

We consider a finite-dimensional vector space $V$ with a nondegenerate
symmetric bilinear form $( ~ | ~ )$.
An isotropic subspace $W \subset V$ is called maximal
if it is not strictly contained in another isotropic subspace of $V$.
Using Witt's theorem and its corollaries \cite{Lang} it is easy
to prove that, if $V$ is a real vector space of signature $(n,n)$
or a complex vector space of dimension $2n$,
the dimension of any maximal isotropic subspace, $W$, of
$V$ is equal to $n$, and that
an isotropic subspace is maximal if and only if it is equal
to its own orthogonal. 
One can also show that
any maximal isotropic subspace $W\subset V$  has isotropic
                                                                              complements, and
that they are maximal. A choice of such a maximal isotropic 
complement, $W' \subset V, V= W \oplus W'$, determines
an isomorphism between the space $W^*$ dual to $W$ and the space $W'$.

If, in addition, we assume that the spaces $V$ and $W$ possess
a Lie-algebra structure,
we arrive at the definition of a Manin pair.

\begin{definition}
A {\it Manin pair} is a pair, $(\d, \g)$, where $\d$ is a Lie algebra   
of even dimension $2n$, with an invariant,  
nondegenerate symmetric bilinear form, 
of signature $(n,n)$ in the real case,
and $\g$ is both
a maximal isotropic subspace and a Lie subalgebra of $\d$.
\end{definition}

\begin{definition}
A {\it Manin quasi-triple} $(\d, \g, \h)$ is a Manin pair $(\d, \g)$
with an isotropic complement $\h$ of $\g$ in $\d$.
\end{definition}

Thus, in a Manin quasi-triple, 
$\d= \g \oplus \h$, $\g$ is a maximal isotropic Lie subalgebra and 
$\h$ is a maximal isotropic linear subspace
of $\d$.
We shall denote by $f: \bigwedge^2 \g
\rightarrow \g$ the Lie-algebra structure of $\g$ which is the restriction of
that of $\d$.  We shall denote by $j: \g^* \rightarrow \h$ the 
isomorphism of vector spaces defined by the decomposition $\d= \g \oplus
\h$, which satisfies
\begin{equation} \nonumber
(j(\xi)|x) = ~ < \xi, x > \ ,
\end{equation}
for each $\xi \in \g^{*}$, and $x \in \g $,
and we shall denote
the projections from $\d$ to $\g$ and $\h$ by $p_{\g}$ and $p_{\h}$
respectively. 
On $\g$, we introduce the cobracket, $F:\g \rightarrow \bigwedge^2
\g$. It is the transpose of the map from $\bigwedge^{2} \g^{*}$ to
$\g^{*}$, which we denote by the same letter, defined by
\begin{equation} \nonumber
F(\xi, \eta) = j^{-1} p_{\h} [j(\xi), j(\eta)], 
\end{equation}
for $\xi, \eta \in \g^{*}$ .
We also introduce an element $\varphi \in \bigwedge^3 \g$
which is defined by the map from $\bigwedge^{2}\g^{*}$ to $\g$, 
denoted again by the same letter,
\begin{equation} \nonumber
\varphi(\xi, \eta)=  p_{\g} [j(\xi), j(\eta)].
\end{equation}
The Lie algebra $\g$, with cobracket $F$ and element $\varphi \in
\bigwedge ^{3} \g$, is called
a {\it Lie quasi-bialgebra} \cite{Dr2}. 
(The element $\varphi$ is the classical limit of
the associator of a quasi-Hopf algebra.)
Conversely, from a Lie quasi-bialgebra $(\g, F, \varphi)$, we obtain
a Manin quasi-triple $(\d, \g , \g^{*})$, where $\d = \g \oplus \g ^{*}$
with its canonical scalar product induced by the pairing between $\g$
and $\g^{*}$, 
\begin{equation} \nonumber
( (x,\xi) | (y,\eta))= < x,\eta > + < y, \xi > ,
\end{equation}
the Lie bracket on $\d$ being defined as follows \cite{Dr2},
\begin{equation} \label{comd}
[x, y] = f(x,y), ~~ \  [x, \xi]= \ad^*_{x} \xi - \ad^*_{\xi} x, ~~ \
[\xi, \eta]= F(\xi, \eta) + \varphi (\xi, \eta),
\end{equation}
where $x,y \in \g$ and $\xi, \eta \in \g^{*}$.
Sometimes we refer to the bialgebra data $F$ and $\varphi$ corresponding
to the complement $\h \subset \d$ as $F_{\h}$ and $\varphi_{\h}$.

\begin{example} 
For any Lie algebra $\g$, the choice $F=0$ and $\varphi=0$ defines the
Manin pair $(\d, \g)$, where $\d=\g \oplus \g^*$ with the Lie
bracket
\begin{equation} \nonumber
[x, y] = f(x,y), ~~ \  [x, \xi]= \ad^*_{x} \xi,~~ \ [\xi, \eta]= 0.
\end{equation}
We call this Manin pair the {\it standard Manin pair} associated to $\g$.
In the standard Manin pair, $\g$ possesses a canonical complement,
$\h=\g^*$,
defining a Manin triple $(\d, \g, \g^*)$.
\end{example}

\begin{example} \label{gminus}
If a Lie algebra $\g$ possesses an invariant, nondegenerate symmetric
bilinear form $K$, one can construct another Manin pair with
$\d= \g \oplus \g$, the direct sum of two copies of 
Lie algebra $\g$. (See \cite{STS1}.) 
The scalar product on $\d$ is defined as the
difference of the bilinear forms on the two copies of $\g$,
\begin{equation} \label{scal}
((x_1,x_{2}) | (y_{1},y_2))= K(x_1,y_{1}) - K(x_{2},y_2) \ ,
\end{equation}
and $\g$ is embedded into $\d$ by the diagonal embedding
$\Delta: x \mapsto (x,x)$. A possible choice of an isotropic
complement to $\g$ is given by $\g_- ={\frac 12}  \Delta_-(\g)$, where
$\Delta_-: x \mapsto (x,-x)$ is the anti-diagonal embedding.
In general, the isotropic subspace $\g_- \subset \d$ is not a
Lie subalgebra, and $(\d, \g, \g_-)$ is only a Manin quasi-triple.
Map $F$ vanishes because $[\g_{-}, \g_-] \subset \g$.
If we identify elements in $\g$ with elements in $\g_{-}$ by
${\frac 12} \Delta_{-}$, and
then with elements in $\g^{*}$ by the isomorphism $j^{-1}$ of this Manin
quasi-triple, the element $\varphi \in \bigwedge ^{3} \g$ can be
identified with the trilinear  form on $\g$,
\begin{equation} \nonumber
(x,y,z) \mapsto {\frac 14} K(z,[x,y]),
\end{equation}
which is actually anti-symmetric and $\ad(\g)$-invariant. Let
$(e_i), i=1, \dots, n$, be a basis of $\g$, and let $(K^{ij})$
be the matrix inverse to matrix $(K_{ij})$, where 
$ K_{ij}= K(e_{i}, e_{j})$. Let us
denote the structure constants of $\g$ in this basis by
$f_{ij}^{k}$. Then
\begin{equation} \nonumber
\varphi = 
\frac{1}{4} K^{il}K^{jm}f_{lm}^{k} e_{i} \otimes e_{j} \otimes e_{k} 
= \frac{1}{24} K^{il}K^{jm}f_{lm}^{k} e_{i} \wedge e_{j} \wedge e_{k}. 
\end{equation}
If, in particular, we choose an orthogonal basis such that 
$K(e_i,e_j)= \frac{1}{2} \delta_{ij}$, then 
\begin{equation} \nonumber
\varphi = \sum_{ijk}
 f_{jk}^{i} e_i\otimes e_j \otimes e_k
= \frac{1}{6} \sum_{ijk}
 f_{jk}^{i} e_i\wedge e_j \wedge e_k.
\end{equation}
\end{example}

\subsection{Twisting}
All Manin quasi-triples  corresponding to the same
Manin pair $(\d, \g)$ differ by a twist $t \in \bigwedge^2 \g$, defined
as follows. An isotropic complement to $\g$ in $\d$ always exists and 
is by no means unique. 
Let $\h$ and $\h'$ be two isotropic complements of $\g$ in $\d$,
and let $1_{\d}=p_{\g} + p_{\h}$ and $1_{\d}=p'_{\g} + p'_{\h'}$ be the
decompositions of the identity map of $\d$ into the sum of the projections
defined by the direct decompositions $\d=\g \oplus \h$ and $\d=\g \oplus \h'$.
Then $1_{\h}= \kappa + \lambda$, where $\kappa$ and $\lambda$ are
the restrictions to $\h$ of projections $p'_{\g}$ and $p'_{\h'}$. 

Let $j: \g^* \rightarrow \h$ and $j': \g^* \rightarrow \h'$ be
the isomorphisms of vector spaces defined by these direct decompositions.
We consider the linear map from $\g^*$ to $\d$, 
called the {\it twist} from $\h$ to $\h'$, 
\begin{equation} \label{twist}
t= j'-j \ .
\end{equation}
We first remark
that $t$ takes values in $\g$. In fact, $j'= \lambda \circ  j$ 
and $t=-\kappa \circ j$. 
It is easy to show that $t$ is anti-symmetric. In fact, for $\xi, \eta
\in \g^{*}$,
\begin{equation} \nonumber
< t(\xi), \eta > + < \xi, t(\eta) > ~ = 
(t(\xi) | j(\eta)) + (j'(\xi) | t(\eta)) = (j(\xi) | j(\eta))= 0,
\end{equation}
where we have used the isotropy of $\g, \h$ and $\h'$.
Hence, the map $t$ defines an element in $\bigwedge^2 \g$ which we denote by
the same letter. By convention, $t(\xi, \eta) = < t(\xi) , \eta >$.

\subsection{The canonical $r$-matrix}

Let $(\d, \g, \h)$ be a Manin quasi-triple. We identify $\d$ with $\g \oplus
\g^*$ using the isomorphism $j^{-1}$ of $\h$ onto $\g^*$. The map
$r_{\d}: \d^* \rightarrow \d$ defined by $r_{\d}: ( \xi , x ) \mapsto (0, \xi)$
for $x \in \g, \xi \in \g^*$ defines an element $r_{\d} \in \d \otimes \d$,
called the {\it canonical} $r$-{\it matrix}.
Let $( e_i), i=1, \dots, n,$ be
a basis of $\g$ and $( \e^i ), i=1, \dots, n,$ be 
the dual basis in $\g ^{*}$. Then
\begin{equation}   \label{rd}
r_{\d}= \sum_{i=1}^n e_i \otimes \e^i.
\end{equation}

The symmetric part of $r_{\d}$ coincides with the
scalar product of $\d$ up to a factor of $2$,
and therefore it is $\ad(\d)$-invariant.
The element $r_{\d}$ satisfies a relation that generalizes
the classical Yang-Baxter equation.
To derive that relation we introduce the notion of a Drinfeld bracket
for elements of $\mathfrak a \otimes \mathfrak a$, for any Lie algebra
$\mathfrak a$.

\begin{definition}
Let $\mathfrak a$ be a Lie algebra with a basis $( e_{\alpha} ), \alpha
= 1 , \ldots , N $, and let 
$r= \sum_{\alpha \beta} r^{\alpha \beta} e_{\alpha} \otimes e_{\beta}$ 
be an element of $\mathfrak a \otimes \mathfrak a $.
The {\it Drinfeld bracket} of $r$ is the element in $\mathfrak a \otimes
\mathfrak a \otimes \mathfrak a$
defined as follows,
\begin{equation} \label{YBb}
\l r, r \r= [r_{12}, r_{13}] + [r_{12}, r_{23}] + [r_{13}, r_{23}],
\end{equation}
where 
$r_{12} = r \otimes 1, 
r_{13} = \sum_{\alpha \beta} r^{\alpha \beta} e_{\alpha} \otimes 1 \otimes
e_{\beta}, 
r_{23} = 1 \otimes r$, and $1$ is the unit of the universal enveloping
algebra $U(\mathfrak a)$. 
\end {definition}

If $r$ is anti-symmetric, then
\begin{equation} \label{skew}
\l r, r\r = - \frac{1}{2} [r,r] \ ,
\end{equation}
where $[ ~, ~ ]$ is the algebraic Schouten bracket 
of $\bigwedge \mathfrak a$. (See, {\it e.g.}, \cite{YKS2} or \cite{YKS3}.)
If $r$ is symmetric and $\ad(\mathfrak a)$-invariant, then
\begin{equation}
\l r, r \r = [r_{12}, r_{13}] \ ,
\end{equation}
and $\l r, r \r$ is the $\ad(\mathfrak a)$-invariant element in
$\bigwedge ^{3} \mathfrak a$ with components 
$r^{\kappa \beta} r^{\lambda \gamma} f^{\alpha}_{\kappa \lambda}$.
For any $r \in \mathfrak a \otimes \mathfrak a$ with 
$\ad(\mathfrak a)$-invariant
symmetric part $s$,
\begin{equation} \label{ras}
\l r, r \r = \l a, a \r + \l s , s \r \ ,
\end{equation}
where $a$ is the anti-symmetric part of $r$. 
Therefore, when the symmetric part of $r$ is  
invariant, $\l r, r \r$ is in $\bigwedge ^{3} \mathfrak a$.

\begin{proposition} \cite{BKS}
The canonical $r$-matrix, $r_{\d}$, associated to the Manin quasi-triple
$(\d, \g, \h)$ satisfies
\begin{equation} \label{cYB}
\l r_{\d}, r_{\d} \r = \varphi,
\end{equation}
where the element $\varphi \in \bigwedge^3 \g$ 
is considered as an element in $\bigwedge ^{3} \d$ \ .
\end{proposition}

\begin{proof}
In the basis $(e_i, \e^i)$, the commutation relations of $\d$ have the form,
\begin{equation} \nonumber
[e_i, e_j]=f_{ij}^k e_k \ , 
\ [e_i, \e^j]=- f_{ik}^j \e^k + F_i^{jk} e_k \ ,
\ [\e^i, \e^j]= F^{ij}_k \e^k + \varphi^{ijk} e_k \ .
\end{equation}
The computation of the three terms entering the Drinfeld bracket yields
\begin{gather} \nonumber
[r_{12}, r_{13}]= [e_i \otimes \e^i \otimes 1, e_j \otimes 1 \otimes \e^j]=
f_{ij}^k e_k \otimes \e^i \otimes \e^j ;
\\
\nonumber 
[r_{12}, r_{23}] =[e_i \otimes \e^i \otimes 1, 1 \otimes e_j \otimes \e^j]=
f_{jk}^i e_i \otimes \e^k \otimes \e^j - 
F_j^{ik} e_i \otimes e_k \otimes \e^j ;
\\
\nonumber
[r_{13}, r_{23}]= [e_i \otimes 1 \otimes \e^i, 1 \otimes e_j \otimes \e^j]=
F_k^{ij} e_i \otimes e_j \otimes \e^k + 
\varphi^{ijk} e_i \otimes e_j \otimes e_k,
\end{gather}
whence 
$\l r_{\d}, r_{\d} \r =\varphi^{ijk} e_i \otimes e_j \otimes e_k=\varphi$.
\end{proof}

We now list properties of the
anti-symmetric part of the canonical $r$-matrix, $r_{\d}$.

\begin{proposition} \label{algSch}
Let $(\d,\g,\h)$ be a Manin quasi-triple, where we identify
$\h$ with $\g^*$. Let $a_{\d} \in  \bigwedge ^{2} \d$ be the anti-symmetric
part of $r_{\d}$. For any $x\in \g, ~ \xi \in \g^*$,
\begin{equation} \label{ad}
[x, a_{\d}] = F(x), ~~ \ [\xi, a_{\d}]= - f(\xi) + \varphi(\xi), ~~ \ 
[a_{\d}, \varphi]=0.
\end{equation}
\end{proposition}

\begin{proof}
The proof of the first two equalities is by computation, 
using the formula
$a_{\d} = \frac{1}{2} \sum _{i =1}^{n} e_i \wedge \e^i$ 
and the derivation property
of the algebraic Schouten bracket. To prove the third equality, we use
formulas \eqref{ras}, \eqref{cYB} and \eqref{skew} to obtain
\begin{equation} \nonumber
[a_{\d}, \varphi] = - \frac{1}{2} [a_{\d},[a_{\d},a_{\d}]] +[a_{\d}, \l
s_{\d}, s_{\d} \r].
\end{equation}
The first term vanishes by the graded Jacobi identity and the second
term vanishes because $\l s_{d}, s_{\d} \r$ is $\ad (\d)$-invariant.
\end{proof}

Under a twist $t$, the canonical $r$-matrix $r_{\d}$  is modified to
\begin{equation} \label{twistr}
r'_{\d}=r_{\d} +t.
\end{equation}
Indeed, any element of $\d$ can be decomposed in two ways, as 
$ x + j(\xi)$ and as $ x' + j'(\xi')$, where 
$x,x' \in~\g, ~ \xi , \xi' \in \g^*$.
Then, by definition, $t(\xi) = j'(\xi ') - j(\xi) \in \g$ , while 
\begin{equation}  \nonumber
(r'_{\d} - r_{\d}) ( j(\xi)+x) =
r'_{\d} (j'( \xi')+x') - r_{\d} (j( \xi)+x)= j'(\xi ') - j(\xi) = t(\xi).
\end{equation}
If $t= \sum_{ij} t^{ij} e_i \otimes e_j
= \frac{1}{2} \sum_{ij} t^{ij} e_i \wedge e_j$, 
then after twisting by $t$, the dual basis of $\g^{*}$ becomes
\begin{equation}
{\e '} ^{i} = \e ^{i} + t^{ij} e_{j} \ .
\end{equation}
Under the  twist $t$ relating isotropic complements $\h$ and $\h'$, 
the cobracket $F$ of the Lie quasi-bialgebra and the element $\varphi$ 
are modified as follows  \cite{Dr2} \cite{YKS1}:
\begin{gather}  \label{F'}
F_{\h'}= F_{\h} + F_{1} ,  \\
\label{phi'}
\varphi_{\h'}= \varphi_{\h} - \l t,t \r + \varphi _{1} \ , 
\end{gather}
where $F_{1}(x)=  \ad_{x}t$ and 
$\varphi _{1} (\xi) =  {\overline{\ad}} ^{*}_{\xi} t$.
Here $\ad$ denotes the adjoint action of $\g$ on $\bigwedge ^{2} \g$, 
while ${\overline{\ad}} ^{*}_{\xi} t$ denotes the projection  
of $\ad^{*}_{\xi}t$ onto $\bigwedge ^{2} \g$, where $\g^{*} \subset \d^{*}$
acts on $\bigwedge ^{2} \g \subset \bigwedge ^{2} \d$ by the coadjoint
action. In fact,
\begin{equation}
\varphi_{1}^{ijk} = F^{jk}_{l} t^{il} - F^{ik}_{l} t^{jl} \ .
\end{equation}

\section{Group pairs and quasi-Poisson structures  }

In this section, we study the global objects corresponding to Manin
pairs and Manin quasi-triples.
In the rest of this paper, we shall assume that the base field is $\R$, 
and all manifolds and maps will be assumed to be smooth.
We shall often abbreviate ``fields of multivectors'' to ``multivectors''.
A {\it bivector map} between $(M_{1}, P_{M_{1}})$ and $(M_{2},
P_{M_{2}})$,
where $P_{M_{i}}$ is a bivector on manifold $M_{i} ~ (i = 1,2)$,
is a map $u$ from $M_{1}$ to $M_{2}$ such that
\begin{equation} \nonumber
u_{*}P_{M_{1}} = P_{M_{2}} \ .
\end{equation}
(If $P_{M_{1}}
$ and $P_{M_{2}}$ are Poisson bivectors, such a map is
called a {\it Poisson map}.)
 
\subsection{Group pairs and quasi-triples}
We first introduce the objects which integrate Manin pairs and Manin
quasi-triples.

\begin{definition} A {\it group pair} is a pair $(D,G)$, where
$D$ is a connected  Lie group with a bi-invariant scalar product
and $G$ is a connected, closed Lie subgroup of $D$, such that the Lie algebras,
$\d$ and $\g$,
of $D$ and $G$ form a Manin pair.
\end{definition}

It is evident from this definition that, given a 
finite-dimensional Manin pair $(\d,\g)$ it can
be integrated into a unique group pair, where $D$ is simply
connected, provided $\g$ is the Lie algebra of a closed Lie 
subgroup of $D$. 

\begin{example}
A group pair corresponding to the standard Manin
pair $(\g \oplus \g ^{*},\g)$ of the Lie algebra $\g$ is $(T^*G,G)$.
Here $T^*G$ is the cotangent bundle of the connected, simply
connected Lie group corresponding to the Lie algebra $\g$, equipped
with the group structure of a semi-direct product,
upon identification with $G \times \g^*$ by left translations.
The group $G$ is embedded into $T^*G$ as the zero section.
\end{example}

\begin{example} \label{diag}
For a Lie algebra $\g$ with an invariant, nondegenerate symmetric bilinear
form, the Manin pair $(\g \oplus \g, \g)$ of Example \ref{gminus} 
can be integrated
to a group pair $(G\times G, G)$. The group $G$ is embedded into
$D=G\times G$ as the diagonal. 
\end{example}

\begin{definition}
A {\it quasi-triple} $(D,G,\h)$ is a group pair $(D,G)$, where an
isotropic complement $\h$ of $\g$ in $\d$ has been chosen.
\end{definition}

Note that a quasi-triple is not a triple of Lie
groups. In general, the subspace $\h \subset \d$ is not a Lie
subalgebra and cannot be integrated into a Lie subgroup.
If $(\d, \g, \h)$ is a Manin triple, integrating $\h$
to a Lie group $H$ yields
a triple of groups $(D,G,H)$. (See \cite{LW}.)

We shall now study the bivectors on $D$ and $G$ which generalize
the multiplicative bivectors obtained when $(\d,\g,\h)$ is a
Manin triple, and $D$ and $G$ are Poisson Lie groups \cite{Dr1}
\cite{STS1}.

\subsection{The quasi-Poisson structure of $D$}
Let $(D,G,\h)$ be a quasi-triple.
We have recalled in Section 2 the definition of the canonical
$r$-matrix of $\d$, $r_{\d} \in \d \otimes \d$. Thus $r_{\d}$
defines a contravariant 2-tensor on $D$,
\begin{equation}  \label{PD}
P_D= r_{\d}^{\lambda} - r_{\d}^{\rho},
\end{equation}
where $r_{\d}^{\lambda}$ (resp., $r_{\d}^{\rho}$) denotes the
left- (resp., right-) invariant 2-tensor on the
Lie group $D$ with value $r_{\d}$
at the identity. Sometimes we denote the bivector $P_D$ corresponding
to the complement $\h \subset \d$ by $P_D^{\h}$ to render the
dependence on $\h$ explicit.

Actually, $P_D$ is a bivector because the symmetric part of $r_{\d}$, 
being $\ad(\d)$-invariant, is invariant 
with respect to the adjoint action of $D$ and
cancels in \eqref{PD}. 

At the identity of $D$, $P_D$ vanishes because there $r_{\d}^{\lambda}$
coincides with $r_{\d}^{\rho}$.  The following multiplicativity
property of $P_D$ is obvious from the definition.

\begin{proposition}
The bivector $P_D$ on the Lie group $D$ is multiplicative
with respect to the multiplication of $D$, {\it i.e.},
the multiplication map $m: D\times D \rightarrow D$ is a bivector map
from $D \times D$ with the
product bivector to $(D,P_D)$.
\end{proposition}

In general, the Schouten bracket $[P_D, P_D]$ does not vanish.

\begin{proposition} \label{propSchPD} 
The Schouten bracket of the bivector $P_{D}$ is given by
\begin{equation} \label{SchPD}
\frac{1}{2} [P_D, P_D] = \varphi^{\rho} - \varphi^{\lambda}.
\end{equation}
\end{proposition}

\begin{proof}
Since $P_{D} = a_{\d}^{\lambda} - a_{\d}^{\rho}$, 
and since left- and right-invariant vector 
fields commute with each other, we obtain
\begin{equation} \nonumber
[P_D, P_D] 
= [a_{\d}^{\lambda} - a_{\d}^{\rho} ,  a_{\d}^{\lambda} - a_{\d}^{\rho}] 
=[a_{\d}^{\lambda}, a_{\d}^{\lambda}] + [a_{\d}^{\rho}, a_{\d}^{\rho}]
= [a_{\d} , a_{\d} ]^{\lambda} - [a_{\d} , a_{\d}]^{\rho} \ .
\end{equation}
Using \eqref{skew}, \eqref{ras} and the fact that $\l s_{\d}, s_{\d} \r$
is $\ad(\d)$-invariant, the right-hand side is $ 2 ( \l r_{\d},
r_{\d} \r^{\rho} - \l r_{\d} , r_{\d} \r ^{\lambda} )$. 
Since $\l r_{\d} ,
r_{\d}\r = \varphi$, we arrive at formula \eqref{SchPD}.
\end{proof}
 
We collect some properties of the bivector $P_D$ in the following
proposition.

\begin{proposition} \label{propPD}
The bivector $P_D$ satisfies
\begin{equation}
\L_{x^{\lambda}} P_D = F(x)^{\lambda}, \
\L_{x^{\rho}} P_D =  F(x)^{\rho} \ ,
\end{equation}
for $x\in \g$, and
\begin{equation}
\L_{\xi^{\lambda}} P_D = - f(\xi)^{\lambda} + \varphi(\xi)^{\lambda}, \
\L_{\xi^{\rho}} P_D = - f(\xi)^{\rho} + \varphi(\xi)^{\rho} \ ,
\end{equation}
for $\xi \in \h$, and
\begin{equation}
[P_D, \varphi^{\lambda}]=0, \ [P_D, \varphi^{\rho}]=0.
\end{equation}
\end{proposition}

\begin{proof}
Here and below $\h$ is identified with $\g^{*}$.
By Proposition \ref{algSch}, 
\begin{gather} \nonumber
[x^{\lambda}, a_{\d}^{\lambda}] = [x, a_{\d}]^{\lambda} = F(x)^{\lambda},
\\ \nonumber
[\xi^{\lambda}, a_{\d}^{\lambda}] = [\xi, a_{\d}]^{\lambda} = 
- f(\xi)^{\lambda} + \varphi(\xi)^{\lambda},
\\ \nonumber
[a_{\d}^{\lambda}, \varphi^{\lambda}]=[a_{\d}, \varphi]^{\lambda}=0.
\end{gather}
To conclude we use the fact that $P_{D} = a_{\d}^{\lambda} - a_{\d}^{\rho}$
and, again, the fact that left- and right-invariant vector 
fields commute. The proof for right-invariant vector fields is similar. 
\end{proof}

Modifying the chosen complement $\h$ of $\g$ by a twist $t\in
\bigwedge^2 \g$ leads to modifying the bivector $P_D$
in the following simple way,
\begin{equation}
P_D^{\h'} = P_D^{\h} + t^{\lambda} - t^{\rho} \ ,
\end{equation}
since $r_{\d}$ is modified according to \eqref{twistr}.

\subsection{The quasi-Poisson structure of $G$}

The bivector $P_D$ has a natural restriction to the subgroup $G
\subset D$. An isotropic complement $\h$ 
such that $\d = \g \oplus \h$ being chosen, to any $g \in G$,
one can associate another 
splitting of $\d$, namely $\d = \g \oplus \Ad_{g} \h$. Here we
use the fact that the Lie subalgebra $\g$ is stable with respect
to $\Ad_{g}$. 
The identification $j'$ of $\g^{*}$ with $\h' = \Ad _{g} \h$ is 
$j' = \Ad_{g} ~ \circ ~ j ~ \circ \ ^t{\Ad _{g }}$, and therefore 
$r_{\d}$ is modified to 
\begin{equation} \label{newr}
r'_{\d} = \Ad_{g} r_{\d} \ ,
\end{equation}
where $\Ad$ denotes the adjoint action of $D$ on the tensor
product $\d \otimes \d$. 
We denote the corresponding twist by $t_g$. 
According to \eqref{twistr} we obtain
\begin{equation}  \label{jg}
t_g =  \Ad_{g} r_{\d} -r_{\d} \ .
\end{equation}
For each $g \in G$, $t_g$ is an element 
of $\bigwedge^2 \g$. Thus, we can define a bivector $P_{G}$ 
on $G$ by right translation of $t_g$, and then
\begin{equation}
P_G= r_{\d}^{\lambda} - r_{\d}^{\rho}.
\end{equation}
Therefore, the embedding of $(G,P_G)$ into $(D,P_D)$ 
is a bivector map. Actually, this requirement
determines $P_G$ uniquely. 

It is clear that $P_G$ inherits the multiplicativity property
of $P_D$ and the Schouten bracket of $P_G$ is given by the same
formula,
\begin{equation}  \label{SchPG}
\frac{1}{2}[P_G , P_G]= \varphi^{\rho} - \varphi^{\lambda}.
\end{equation}
Note that both $\varphi^{\rho}_g$ and $ \varphi^{\lambda}_g$ are already
elements of $\bigwedge^3 T_gG$, whereas for a canonical $r$-matrix
it is only the difference $ (r_{\d}^{\lambda}-
r_{\d}^{\rho})_g$ which lies in $\bigwedge^2 T_gG$.
Moreover,
\begin{equation} \label{pentagon}
[P_{G}, \varphi^{\lambda}] = 0, ~~ \ [P_{G}, \varphi^{\rho}] = 0 \ .
\end{equation}

We now recall from \cite{YKS1} and \cite{YKS2} 
the definition of quasi-Poisson Lie groups, 
which we shall use in the study of
quasi-Poisson actions in Section 4.

\begin{definition}
A {\it quasi-Poisson structure} on a Lie group $G$ 
is defined by a multiplicative
bivector  $P_G$ and an element $\varphi$ in $\bigwedge^3 \g$ such that
$\frac{1}{2} [P_G,P_G]=\varphi^{\rho}-\varphi^{\lambda}$ and 
$[P_G,\varphi^{\lambda}]= 0.$
\end{definition}

It also follows from the definition that 
$[ P_{G} , \varphi ^{\rho} ] = 0$, because, by 
the graded Jacobi identity, 
$[ P_{G} , \varphi^{\rho} - \varphi ^{\lambda} ] = 0$.
Identity \eqref{pentagon}
is the classical limit of the
pentagon identity for quasi-Hopf algebras \cite{Dr2}.
For a connected Lie group $G$, the multiplicativity condition
is equivalent to its infinitesimal version,
\begin{equation}
\L_{x^{\lambda}} P_G=  F(x)^{\lambda},
\end{equation}
for each $x\in \g$.

\medskip

Equations \eqref{SchPG} and \eqref{pentagon} show that
the Lie group $G$, equipped with the bivector $P_G$ and
the element $\varphi$ in $\bigwedge^3 \g$, is a quasi-Poisson Lie group,
in the sense of the preceding definition, and so is the Lie group $D$,
with the bivector $P_{D}$ and the same element $\varphi$ in $\bigwedge^3
\g$, considered as an element in $\bigwedge^{3} \d$. This last fact 
follows from the properties
of $P_D$ stated in Propositions \ref{propSchPD} and \ref{propPD}.

We sometimes refer to the bivector $P_{G}$ as $P_{G}^{\h}$ to emphasize the
dependence on the choice of a complement $\h$. 
For the bivector $P_G^{\h}$ corresponding to the complement $\h$,
all the properties
required in the definition of a quasi-Poisson Lie group follow from the
analogous properties of $P^{\h}_D$, considered in the previous section.
Hence, $(G,P^{\h}_G,\varphi_{\h})$ is a quasi-Poisson Lie group
integrating the Lie quasi-bialgebra structure of $\g$, defined by the
Manin quasi-triple $(\d,\g,\h)$. We often refer to 
$(G,P^{\h}_G,\varphi_{\h})$ as $G_D^{\h}$.

Bivector $P_G^{\h}$ vanishes if the complement $\h \subset \d$ is
$\ad(\g)$-invariant, $[\g, \h]\subset \h$, because from
$\Ad_g \h =\h$ it follows that $t_g = \Ad_{g} r_{\d} - r_{\d}=0$.
This simple observation implies that, 
in the standard quasi-triple $(T^*G,G,\g^*)$,
the bivector $P_G = P_{G}^{\g^{*}}$ on the group $G$ vanishes,
and also in the quasi-triple $(G\times G, G, \g_-)$, the bivector 
$P_{G} = P_{G}^{\g_{-}}$ vanishes.

Under a twist $t$, the bivector $P_G^{\h}$ is modified in the same way as
the bivector $P_D^{\h}$,
\begin{equation}  \label{P'G}
P_G^{\h'} = P_G^{\h} + t^{\lambda} - t^{\rho}.
\end{equation}

\subsection{The dressing action of $D$ on $D/G$}

To any group pair $(D,G)$ one can associate the quotient
space $D/G$, which we shall denote by $S$. The space $S$ 
will be the target of the generalized moment maps
in Section 5.

The action of $D$ on itself by left multiplication 
induces an action of $D$ on $S$. 
Because this action generalizes the dressing of group-valued
solutions of field equations, it is called the `dressing action'
\cite{STS1} \cite{STS2}. We denote the corresponding infinitesimal action by 
$X \rightarrow X_{S} $, for $X \in \d$. By definition, 
$X_{S}$ is the projection onto $S$ of the 
opposite of the right-invariant vector field on $D$
with value $X$ at the identity, $e$.

\begin{definition}
An isotropic complement $\h$ to $\g$ in $\d$ is called {\it admissible}
at the point $s\in D/G$ if the infinitesimal dressing 
action restricted to $\h$ defines an isomorphism from
$\h$ onto $T_s(D/G)$.
\end{definition}

It is clear that any isotropic complement $\h$ to $\g$ is
admissible in some open neighborhood of $eG \in D/G$.
If the complement $\h$ is admissible at a point $s\in D/G$,
it is also admissible in some open neighborhood $U$ of $s$.
If there exists an $\h$ which is admissible everywhere on $D/G$, we call
the corresponding quasi-triple $(D,G,\h)$ {\it complete}.

If we identitfy the tangent spaces to $D$ with $\d $ by means of right
translations, then the tangent space to $D/G$ at $s \in D/G$ is
identified with $\d / \Ad_{s}\g$, and, for $X \in \d$, $X_{S}(s)$ 
is identified with the
class of $X$ in $\d / \Ad_{s}\g$.
If $\h$ is admissible, then $\h$ is isomorphic to $T_{s}(D/G)$, and we
obtain both decompositions,
\begin{equation} \nonumber
\d = \g \oplus \h = \Ad_{s} \g \oplus \h \ .
\end{equation}
If $x$ is in $\g$, then $x_{S}(s)$ is identified with the element
$\theta_{s}(x) = (p_{\Ad_{s}\g} - p_{\g} )(x)\in \h$, the difference of
the projections of $x$ onto $\Ad_{s}\g$ and $\g$, parallel to $\h$.
Thus $x_{S}(s) = (\theta_{s}(x))_{S}(s) $.
So, if this map $\theta_{s}$ from $\g$ to $\h$
is composed with the isomorphism $j^{-1}$ from $\h$ to $\g^{*}$ defined
in Section 2.1, we obtain
$\tau_{s} = j^{-1} \circ \theta_{s}$ from $\g$ to $\g^{*}$ satsifying
\begin{equation} \label{maptau}
x_{S}(s) = ((j \circ \tau_{s})(x))_{S}(s) \ 
\end{equation}
and, in particular,
\begin{equation} \label{tau}
(e_{i})_{S} (s) = (\tau_{s})_{ik}(j\e^{k})_{S}(s)  .
\end{equation}
Because $\g, \h$ and $\Ad_{s}\g$ are isotropic, $\tau_{s}$ is
anti-symmetric and defines an element in $\bigwedge^{2}\g^{*}$, which we
denote by the same letter.

\begin{proposition}
At any point $s \in D/G$
there exists an admissible complement $\h$ of $\g$ in $\d$.
\end{proposition}

\begin{proof}
Maximal isotropic subspaces in $\d$ form a Grassmannian
which we denote by $\G(\d)$. Since $D$ is a connected Lie group, the
subspaces $\g$ and $\Ad_{s}\g$ belong to the same connected component
of $\G(\d)$. Let $\h'$ be an isotropic complement of $\g$ in $\d$. 
The Grassmannian $\G(\d)$ being an algebraic variety, the set of isotropic
complements $\h$ of $\g$ in the connected component of $\h'$ is a
Zariski open set. Since $\Ad_{s} \g$ is in the same connected component
of $\G(\d)$ as $\g$, the set of isotropic complements to
$\Ad_{s} \g$ in the connected component of $\h'$ is also a
Zariski open set.
An intersection of two nonempty Zariski open sets being
nonempty, one can always find an $\h$ which is an isotropic
complement of both $\g$ and $\Ad_{s} \g$. Any such subspace
is admissible at the point $s$.
\end{proof}

The choice of an admissible isotropic complement $\h$ of $\g$ in $\d$
at a point $s$ gives
rise to an additional structure on $D/G$. The space $\h$ being
isomorphic to $T_s(D/G)$, we can define a
map from $\g$ to the $1$-forms on $D/G$ at the point $s$, 
$x \mapsto \hat{x}_{\h}(s)$, as follows, 
\begin{equation} \nonumber
< \hat{x}_{\h}(s), \xi_S(s) > ~ = - (x~|~\xi), 
\end{equation}
for each $\xi  \in \h$,
where in the
left-hand side
we have used the canonical pairing between $1$-forms and vectors. 
The map $\h \rightarrow T_s(D/G)$ being an isomorphism
of linear spaces, so is the map $\g \rightarrow T_s^*(D/G)$.
In other words, the forms $\hat{x}_{\h}(s)$ span 
the cotangent space to $D/G$ at $s$. Thus if $\h$ is admissible in an
open neighborhood $U$ of $s$, we define the map $x \mapsto \hat{x}_{\h}$
from $\g$ to
the $1$-forms on $U \subset D/G$, such that 
\begin{equation} \label{hat}
< \hat{x}_{\h}, \xi_S > ~ = - (x~|~\xi) \ . 
\end{equation}
When the quasi-triple $(D, G, \h)$ is complete,
forms $\hat{x}_{\h}$ are globally defined on $D/G$.
When $(\d, \g, \h)$ is a Manin triple, integrated to Lie groups 
$D, G, H$, then $\xi_{S}$ is identified with the opposite of
the right-invariant vector field on $H$ with value $\xi$
at the identity, and therefore $\hat{x}_{\h}$ is the 
right-invariant $1$-form on $H$ with value $x \in \g \simeq \h^{*}$ 
at the identity.

We now study the effect of a twist $t$ on the map $x \mapsto \hat{x}_{\h}$.
Let us choose two complements $\h$ and $\h'$ admissible at $s\in D/G$.
We first compare the dressing vector fields $\xi_{S}$ and $\xi'_{S}$
on $D/G$ defined by $\xi_{0} \in \g^{*}$ corresponding
to the splittings $\d = \g \oplus \h$ and $\d= \g \oplus \h'$,
respectively, where the twist from 
$\h$ to $\h'$ is $t : \g^{*} \rightarrow \g$.
Then $\xi_{S}= (j\xi_{0})_{S}$, $\xi'_{S} = (j'\xi_{0})_{S}$. It follows
from \eqref{maptau} that $(t\xi_{0})_{S}(s) 
= ((j \circ \tau_{s})(t\xi_{0}))_{S}(s)$, 
and therefore
$\xi'_{S}(s) = (j (\sigma_{s}\xi_{0}))_{S}(s)$, where
$\sigma_{s}= 1_{\g^{*}} + \tau_{s} \circ t$.
Thus
\begin{equation} \nonumber
<x, \xi_{0}> ~ = ~ < \hat{x}_{\h}(s), (j \xi_{0})_{S}(s) > ~ = 
~ < \hat{x}_{\h'}(s), (j (\sigma_{s} \xi_{0}))_{S}(s) >.
\end{equation}
The pairing on the left hand side being nondegenerate,
$\sigma_{s}$ is invertible. Then
\begin{equation} \nonumber
< x , \xi_{0} > ~ = ~ < ~ ^{t}\sigma_{s}^{-1}x ~ , ~ \sigma_{s}\xi_{0} > ~
= ~ < \widehat{(^{t}\sigma_{s}^{-1}x)}_{\h}(s),
(j (\sigma_{s} \xi_{0}))_{S}(s) > \ .
\end{equation}
Therefore,
\begin{equation} \label{twisthat}
\hat{x}_{\h'}(s) = \widehat{(\nu_{s}x)}_{\h}(s) \ ,
\end{equation}
where $\nu_{s} = ~ ^{t}\sigma_{s}^{-1} = (1_{\g} +t \circ
\tau_{s})^{-1}$,
because $t$ and $\tau_{s}$ are anti-symmetric.

\begin{example}
For the standard group pair $(T^*G,G)$, the space $S=T^{*}G/G$ coincides
with $\g^{*}$, the dual space to the Lie algebra $\g$.
Using dual bases
$(e_i)$ and $(\e^i)$ of $\g$ and $\g^*$ one can introduce linear coordinates
$(\xi_i)$ on $\g^*$. The vector fields generating the dressing
action are as follows, 
$\e^i_S = - \partial/ \partial \xi_i$, for the action of $\g^*$,
and $(e_i)_S(\xi) = - \ad^{*}_{e_{i}}\xi = 
f_{ij}^k \xi_k \partial/\partial \xi_j$, for the
action of $\g$.
The quasi-triple $(T^*G,G,\g^*)$ 
is complete, and the $1$-forms corresponding
to the elements $e_i \in \g$ 
correspond to the differentials of the linear coordinates,
$\hat{e}_i=d\xi_i$.
\end{example}

\begin{example} \label{ex}
For a Lie algebra $\g$ with an invariant, nondegenerate symmetric
bilinear form,
the group pair $(G\times G, G)$ of Example \ref{diag} defines the space 
$S=(G\times G)/G \cong G$.
One can express the dressing vector fields as 
$(x,0)_S = - x^{\rho}$
and $(0,x)_S= x^{\lambda}$, where $x^{\rho}$ and 
$x^{\lambda}$ are the  right- and left-invariant 
vector fields on the group $G$
defined by the element $x \in \g$, so that $G$ acts on $S \cong G$ 
by the adjoint action. 
The quasi-triple $(G\times G, G, \g_-)$ is not necessarily complete.
The $1$-forms $\hat{x}$ are defined by the equation
$< \hat{x}, {\frac 12} (y^{\lambda} + y^{\rho}) > ~ =  ((x,x)|(y,-y))$. 
Using the ${\Ad}(G)$-invariance of $K$, we obtain  
$$
< \hat{x}, y^{\lambda} > (g) 
= 2 K(x,(1_{\g}+{\Ad}_g^{-1})^{-1}y)=2K((1_{\g}+{\Ad}_g)^{-1}x,y) \ .
$$
Therefore
\begin{equation} \label{example}
\hat{x}_{g}= 2 (K^{\sharp} \circ (1_{\g}+\Ad_{g})^{-1}
x)^{\lambda}_{g} \ ,
\end{equation}
where $K^{\sharp}(x)(y) = K(x,y)$.
It is clear that the $1$-forms $\hat{x}$ are
well defined if $- 1$ is not an eigenvalue of the operator ${\Ad}_{g}$. 

At the points
$g\in G$ such that ${\rm Ad}_g$ has eigenvalue
$-1$, 
the complement $\g_{-}$, defined by 
$\g_- = \{ {\frac 12}(x, -x) | x \in \g \} $,
is not admissible.
Let us assume that $G$ is a compact simple Lie group. 
Without loss of generality, we can assume that $g$
belongs to a maximal torus $T$. The eigenspace
$V$ corresponding to the eigenvalue $-1$ of $\Ad_{g}$ 
in $\g$ is then even-dimensional, and
it splits into the orthogonal direct sum,
$$
V = \oplus_{\alpha \in \Gamma} V_\alpha,
$$
where each $V_\alpha$ is two-dimensional, and spanned by
the linear combination of the root vectors $e_\alpha,
e_{-\alpha}$,
$$
a_\alpha = \frac{e_\alpha + e_{-\alpha}}{\sqrt{2}}, \
~~ b_\alpha = \frac{e_\alpha - e_{-\alpha}}{i \sqrt{2}} \ ,
$$
the sum being taken 
over a subset  ${\Gamma}$ of the set of positive roots.
If $K(e_{\alpha}, e_{-\alpha})=1$, then $(a_{\alpha},b_{\alpha})$ 
is an orthonormal basis of
$V_{\alpha}$.
Let us twist the complement $\g_{-}$ by  
$t = \frac{\varepsilon}{2} \sum_{\alpha \in \Gamma} 
a_{\alpha} \wedge b_{\alpha}$, 
where $\e $ is a nonzero real
number,. 
After the twist, the new complement is the direct sum of the orthogonal of $V$
and of the subsbace of $\oplus_{\alpha} (V_{\alpha} \oplus V_{\alpha})$
spanned by the vectors,
$$\frac{1}{2}  (a_{\alpha}, -a_{\alpha}) + 
\varepsilon (b_{\alpha}, b_{\alpha}) , \
~~\frac{1}{2} (b_{\alpha}, -b_{\alpha}) 
- \varepsilon (a_{\alpha}, a_{\alpha}) .  
$$
Let us 
compute the $1$-forms ${\hat x}'$ defined by the choice of the 
admissible complement thus obtained by twisting.
If $x$ is in the orthogonal complement to the eigenspace $V$, then 
${\hat x}'$ coincides with the $1$-forms $\hat x$ of formula
\eqref{example}. 
In order to determine the $1$-forms corresponding to a
basis of $V$, it is sufficient to deal with one of the 
$V_{\alpha}$'s,
and we shall omit the index $\alpha$. 
The dressing vector field on $G$ corresponding to 
$(x,x)$  
is $x^{\lambda}-x^{\rho}$,
and that corresponding to $(x,-x)$ is
$-(x^{\lambda}+x^{\rho})$. Assuming that $x$ is in $V$, then 
$\Ad_{g} x = -x$,
and therefore $(x,x)_{G}=2 x^{\lambda}$
and $(x,-x)_{G}=0$. Thus, the dressing vector fields
corresponding to this basis of $V$ are 
$2{\varepsilon} b^{\lambda}$ and    
$-2 {\varepsilon} a^{\lambda}$, respectively.
The dual $1$-forms are defined by 
$$
< \hat a', 2 {\varepsilon} b^{\lambda} > ~ 
= - \left ( (a,a)| (\frac{1}{2} a+\varepsilon
b, -\frac{1}{2}a + {\varepsilon} b) \right)
= - K(a,a) \ ,
$$
$$
< \hat a', 2 {\varepsilon} a^{\lambda} > ~ 
= \left ( (a,a)| (\frac{1}{2} b-\varepsilon
a, -{\frac 12} b - \varepsilon a)\right ) = K(a,b)
$$
and similarly for b,
$$
<\hat b', 2{\varepsilon}b^{\lambda}> ~ = - K(a,b) \ , 
$$
$$
<\hat b', 2  {\varepsilon} a^{\lambda}> ~ = K(b,b) \ .
$$
Therefore 
$$
\hat{a}' = - \frac {1}{2\varepsilon} K(b,\theta), \
~~ \hat{b}' = \frac{1}{2\varepsilon} K(a,\theta),
$$
where $\theta$ is the left-invariant Maurer-Cartan form.
\end{example}

\subsection{The bivector on  $D/G$}

If we choose a quasi-triple $(D,G,\h)$ corresponding to the group pair 
$(D,G)$, a bivector is defined on
the space $S=D/G$ introduced in the previous subsection.
Since the bivector $P_D = r_{\d}^{\lambda} - r_{\d}^{\rho}$ 
is projectable by the projection of $D$ onto $D/G$, 
it defines a bivector $P_{S}$ on $S=D/G$. 
Because all left-invariant 
vector fields generated by $\g$ are projected to zero, the projection of 
$ r_{\d}^{\lambda} $ vanishes, and therefore
\begin{equation}
P_S = - (r_{\d})_{S},
\end{equation}
The notation $(r_{\d})_{S}$ refers to the 
homomorphism from $\d$ to the vector fields on $ D/G$
induced by the action of $D$ on $D/G$, extended to a homomorphism from
the tensor algebra of $\d$ to the algebra of
contravariant tensors on $S=D/G$. 
Note that although
$r_{\d} \in \d \otimes \d $ is not anti-symmetric, after projection
it defines a bivector on $S$.
The bivector $P_S$ 
corresponding to a complement $\h$ will often be referred
to as $P_S^{\h}$.

The properties of $P_S$ follow from those of $P_D$.
Since the Schouten bracket of $P_{S}$ is the projection of that of
$P_{D}$, from \eqref{SchPD} we obtain,
\begin{equation}
\frac{1}{2} [P_S, P_S]= \varphi_{S}.
\end{equation}

Clearly, the action of $D$ on $D/G$ is a bivector map
from $(D,P_{D}) \times (S,P_{S})$ to $(S,P_{S})$. It follows from 
\eqref{ad} that 
\begin{equation}
\L_{x_S} P_S= - F(x)_S \ , \ \L_{\xi_S} P_S=  f(\xi)_S - \varphi(\xi)_S
\ ,
\end{equation}
for all $x \in  \g, \xi \in \h$.

Under a twist $t$ of the Manin quasi-triple $(\d, \g, \h)$ into  
$(\d, \g, \h')$, the bivector $P^{\h}_{S}$
is modified
to $P^{\h'}_{S}= P^{\h}_{S} - t_{S}$.

If the isotropic complement $\h$ of $\g$ 
in the Manin quasi-triple $(\d, \g ,\h)$ is
a Lie subalgebra of $\d$, all the bivectors $P^{\h}_D$, $P^{\h}_G$ and
 $P^{\h}_S$ have
vanishing Schouten brackets, and define Poisson brackets 
satisfying the Jacobi identity on the corresponding
spaces.

The bivector $P_S$ has an interesting characteristic property
which plays the key role in the moment map theory.
Let $(P^{\h}_S)^{\sharp}$ be the map from $1$-forms to vectors defined by 
$<(P^{\h}_S)^{\sharp}\alpha, \beta > ~= P^{\h}_{S}(\alpha,
\beta)$, for any $1$-forms $\alpha, \beta$ on $S$. 
\begin{proposition}
Let $(D,G,\h)$ be a  quasi-triple such that $\h$ is admissible
on an open neighborhood $U$ of $s \in D/G$,  and let $P_S^{\h}$ be the
corresponding bivector on $S=D/G$. Then, for any $x \in \g$, 
\begin{equation}  \label{Smom}
(P_S^{\h})^{\sharp}(\hat{x}_{\h})=x_S \ ,
\end{equation}
holds on $U$, where $\hat{x}_{\h}$ is defined by \eqref{hat}.
This property uniquely characterizes the bivector $P_S^{\h}$ in the
neighborhood $U$.
\end{proposition}

\begin{proof}
We choose an isotropic complement $\h$ of $\g$ in $\d$, 
admissible in an open neighborhood
$U$ of $s \in D/G$.
Let $x$ be an element in $\g$.
We apply the map $(P_{S}^{\h})^{\sharp}$, where   
$P^{\h}_S=- (e_{i})_{S} \otimes (\e^i)_S $,
to the $1$-form $\hat{x}_{\h}$. By definition, we obtain, 
\begin{equation}
(P^{\h}_{S})^{\sharp} (\hat{x}_{\h})=( x | \e^{i})(e_i)_S = x_S \ ,
\end{equation}
which proves formula \eqref{Smom}.
The complement $\h$ being admissible, the $1$-forms $(\hat{e}_i)_{\h}$
form a basis of the cotangent space to $D/G$ at each point in $U$. 
Hence, formula
\eqref{Smom} gives a characterization of the bivector $P^{\h}_S$ on
$U$.
\end{proof}

\begin{example}
The standard Manin quasi-triple $(\d,\g,\g^*)$ is in fact a Manin triple,
so bivector $P_S$ has a vanishing Schouten bracket and defines a
Poisson bracket on $\g^*$. Actually,
the induced bivector on $S = \g^*$
coincides with the Kirillov-Kostant-Souriau
bivector,
\begin{equation} \label{KKS}
P_S= - (r_{\d})_{S} = - (e_i)_{S} \otimes (\e^i)_S= 
\frac{1}{2} \sum_{ijk} f_{ij}^k \xi_k \ \frac{\partial}{\partial \xi_i} \wedge 
\frac{\partial}{\partial \xi_j}.
\end{equation}
\end{example}

\begin{example} \label{exP}
For the quasi-triple $(G\times G,G,\g_-)$, 
the bivector $P_S$ has the form
\begin{equation} \label{lr}
P_S= - (r_{\d})_{S} = - \sum_{i} (\Delta e_i)_{S} \otimes (\Delta_{-} e_i)_S=
 \frac{1}{2} \sum _{i} e_i^{\lambda} \wedge e_i^{\rho}.
\end{equation}
Here we used the fact that $\sum_{i} e_ie_i \in U(\g)$ is a Casimir element,
and therefore $\sum_{i} e_i^{\lambda}\otimes e_i^{\lambda}=\sum_{i} 
e_i^{\rho} \otimes e_i^{\rho}$.
Usually, the Schouten bracket of the bivector \eqref{lr} is nonvanishing.
A notable exception is the group $G=SU(2)$, where $\varphi_S$ vanishes
although $\varphi \neq 0$ \cite{Volkov}. 
\end{example}

\section{Quasi-Poisson actions}

We shall now introduce quasi-Poisson actions in general,
and show that the actions that we have described in
Section 3 are examples of quasi-Poisson actions arising from
Manin pairs.

\subsection{The definition of quasi-Poisson actions}

Let $G$ be a connected Lie group with Lie algebra $\g$, and
let $M$ be a manifold on the which the Lie group $G$ acts.
We shall denote by $x_M$ the vector field  on $M$ 
corresponding to $x\in \g$, and, more generally, the multivector field on
$M$ corresponding to $x \in \bigwedge \g$. By convention, $x_{M}$
satisfies 
\begin{equation} \label{inf}
(x_{M}.f)(m)= {\frac {d}{dt}}f({\exp}(-tx).m)|_{t=0} \ ,
\end{equation}
for $x \in \g, m \in M$ and $f \in C^{\infty}(M)$.
We first recall from \cite{LW} the following
characterization of Poisson actions of connected Poisson Lie groups
on Poisson manifolds.
\begin{proposition}
If the action of a Poisson Lie group $(G,P_G)$ on a Poisson manifold
$(M,P_M)$ is a Poisson action, then, for each $x \in \g$,
\begin{equation} \label{Pact}
\L_{x_M} P_M = - F(x)_M.
\end{equation}
The converse holds if $G$ is connected.
\end{proposition}

Thus, in this case $[P_M,P_M]=0$ (since $P_M$ is Poisson)
and \eqref{Pact} holds. 

\begin{remark} \label{comm}
There is a simple interpretation of the above characterization
of Poisson actions. Equation \eqref{Pact}
is equivalent to the commutativity of the diagram 
$$
\begin{CD}
\bigwedge \g @ > d_F >>  \bigwedge \g \\
@ V ~ VV @ VV V \\ 
\bigwedge{\mathcal X}(M) @ > d_{P_M} >> \bigwedge{\mathcal X}(M)
\end{CD}
$$
where the vertical arrows are induced by the infinitesimal action 
$x \in \g \mapsto x_M \in {\mathcal X}(M)$, the linear space of 
vector fields on $M$. The map $d_F$ is, up to a sign, 
the Chevalley-Eilenberg cohomology operator of the Lie algebra 
$\g^*$ with bracket $F$, with values in the trivial $\g^*$-module.
More precisely, $d_F = ${\large [}$F,~.~ ${\large ]}, 
where {\large [}$~,~${\large ]} is the ``big bracket'' 
on $\bigwedge(\g \oplus \g^*)$.
(See \cite{YKS2}.) 
The map $d_{P_M} = [P_M,~.~]$, where $[~,~]$ is the 
Schouten bracket of multivectors, is the 
Lichnerowicz-Poisson cohomology operator on multivectors on $M$.
In other words, the action 
$x \in \g \mapsto x_M \in {\mathcal X}(M)$ 
is the infinitesimal of a Poisson action 
if and only if it induces a morphism from the complex 
$(\bigwedge \g , d_F)$ 
to the complex $(\bigwedge{\mathcal X}(M), d_{P_M})$.
\end{remark}

The following definition generalizes the
notion of a Poisson action.

\begin{definition}
Let $(G,P_G,\varphi)$ be a connected quasi-Poisson Lie group acting
on a manifold $M$ with a bivector $P_M$. The action of $G$ on $M$
is said to be a {\it quasi-Poisson action} if and only if
\begin{gather} \label{qP1}
\frac{1}{2} [P_M,P_M]= \varphi_M , \\
\label{qP2}
\L_{x_M} P_M= - F(x)_M,
\end{gather}
for each $x\in \g$.
\end{definition}

Let $(D,G,\h)$ be a quasi-triple. We consider $G$ with the
quasi-Poisson structure defined in Section 3.3, and
the space $S=D/G$ with the bivector $P_S$ defined
in Section 3.5. Then the dressing action of $G$ on
$D/G$ is quasi-Poisson. In fact, properties \eqref{qP1}
and \eqref{qP2} were proved in Section 3.5.
The action of $G$ obtained by restriction to any 
$G$-invariant embedded submanifold of 
$D/G$ is also a quasi-Poisson action.
In particular, those $G$-orbits that are embedded submanifolds of $D/G$
are {\it quasi-Poisson spaces}, in the sense that they are manifolds with a
bivector on which a quasi-Poisson Lie group acts by a quasi-Poisson action.

Under a twist $t$, the quasi-Poisson Lie group 
$(G, P^{\h}_{G}, \varphi_{\h} )$ is
modified to $(G, P^{\h'}_{G}, \varphi_{\h'})$, where 
$P^{\h'}_{G}, \varphi_{\h'}$ 
are given by \eqref{P'G} and \eqref{phi'}, 
and $(M,P^{\h}_{M})$ is modified to $(M,P^{\h'}_{M})$,
where 
\begin{equation} \label{twistP}
P^{\h'}_{M}= P^{\h}_{M} - t_{M} \ .
\end{equation}
A simple computation in terms of Schouten brackets shows that after the
twist the action remains quasi-Poisson,
\begin{gather}  \nonumber
\L_{x_M} P^{\h'}_M = \L_{x_M} P^{\h}_M
 - [x_M, t_M]= - (F_{\h}(x) + \ad_{x}t)_M = - F_{\h'}(x)_M, \\
\nonumber
\frac{1}{2}[P^{\h'}_M, P^{\h'}_M]=
 \frac{1}{2} [P^{\h}_M, P^{\h}_M] + \frac{1}{2} [t_{M},t_{M}]
-  [P^{\h}_M, t_M] 
= (\varphi_{\h} - \l t,t \r + \varphi_{1})_M = (\varphi_{\h'})_M,
\end{gather}
where we have used the transformation rules and the notations of
\eqref{F'} and \eqref{phi'}.

The last calculation shows that one can consider a family of quasi-Poisson Lie
groups $G_D^{\h}$ acting on a family of quasi-Poisson spaces $(M,P_M^{\h})$,
where $P_M^{\h'}=P_M^{\h} - t_M$ when
the complements $\h$ and $\h'$ are related
by a twist $t$. We have just shown that when the action of 
$G_D^{\h}$ on $(M,P_M^{\h})$
is quasi-Poisson for an isotropic complement $\h$ it is also
quasi-Poisson for any $\h'$. In the moment map theory which we shall present 
in the next section, it is more convenient to consider families
$G_D^{\h}$ acting on $(M,P_M^{\h})$ than individual quasi-Poisson Lie
groups acting on individual quasi-Poisson spaces.

\begin{remark}
In the case of quasi-Poisson actions, the operators 
$d_F = $ {\large [}$F, ~.~${\large ]} 
and $d_{P_M} = [P_M,~.~]$ can still be 
defined, but their squares do not vanish. 
In fact, $(d_F)^2 =[\varphi,~.~]$ and
$(d_{P_M})^2 = [\varphi_M,~.~]$.
In the first formula the bracket 
is the algebraic Schouten bracket 
on $\bigwedge \g$ and, in the second, it
is the Schouten bracket of multivectors on $M$. 
Formula \eqref{qP2} can still be interpreted as the commutattivity
of the diagram of Remark \ref{comm}, defined by the map
from $\bigwedge \g$ to 
$\bigwedge {\mathcal X}(M)$ induced by the action
of $G$ on $M$, and by the operators
$d_F$ and $d_{P_M}$. It follows that the squares of these operators
commute 
with the induced map, and this implies that the $3$-vector 
${\frac 12} [P,P]- \varphi_M$ has vanishing Schouten bracket with any
multivector in the image of the induced map. Condition 
\eqref{qP1} expresses the fact that this $3$-vector actually vanishes.
\end{remark}

\subsection{Properties of quasi-Poisson actions} \label{actions}
First, we characterize quasi-Poisson actions as bivector maps.

\begin{proposition}
Let $\rho~: G\times M\rightarrow M$ be an action of a connected
quasi-Poisson Lie group $G_D^{\h}$ on a manifold $M$
equipped with a bivector $P^{\h}_M$ which satisfies the
property $\frac{1}{2} [P^{\h}_M,P^{\h}_M]=(\varphi_{\h})_{M}$. 
Then $\rho$ is a quasi-Poisson action
if and only if $\rho$ maps the
bivector $P^{\h}_G+P^{\h}_M$ on $G\times M$ to
$P^{\h}_M$.
\end{proposition}

\begin{proof}
The proof is identical to that in the case of Poisson actions. See
\cite{LW} \cite{YKS3}. 
\end{proof}

Next, we introduce the notion of quasi-Poisson reduction, similar
to the usual Poisson reduction, which yields a genuine 
Poisson structure on the space of orbits.

\begin{theorem} \label{red}
Let $G_D^{\h}$ be a connected quasi-Poisson Lie group acting on
a manifold $(M,P^{\h}_M)$. 
Then the bivector $P_M^{\h}$ defines a Poisson bracket on
the space $C^{\infty}(M)^G$ of smooth $G$-invariant functions
on $M$. This Poisson bracket is independent of the choice 
of $\h$.
\end{theorem}

\begin{proof}
First, we show that the bracket of $G$-invariant
functions, $f_{1}, f_{2}$, is $G$-invariant. Indeed,
\begin{equation} \nonumber
(\L_{x_M} P_M^{\h})(df_1,df_2) = -F(x)_M(df_1,df_2)=0 \ ,
\end{equation}
for any $x\in \g$, because $F(x) \in \bigwedge^2 \g$ and, hence,
$F(x)_{M}$
annihilates $df_1 \wedge df_2$.
Next, we observe that the bracket defined by $P_M^{\h}$ on invariant
functions is a Poisson bracket,
\begin{equation} \nonumber
\frac{1}{2}[P_M^{\h}, P_M^{\h}](df_1,df_2,df_3)= 
\varphi_M(df_1,df_2,df_3)=0 \ ,
\end{equation}
because $\varphi \in \bigwedge^3 \g$ and $f_1, f_2$ and $f_3$ are
$G$-invariant.
Finally, if one modifies the complement $\h$ to $\h'$, the bivector
on $M$ is modified by a twist, $P^{\h'}_M = P_M^{\h} - t_M$, and the 
Poisson bracket of invariant functions is unchanged,
\begin{equation}  \nonumber
P^{\h'}_M(df_1,df_2)=  P_M^{\h}(df_1,df_2)- t_M(df_1,df_2)= P_M^{\h}(df_1,df_2)
\end{equation}
because $t_M$ annihilates invariant functions.
\end{proof}

Let us introduce the projection $p: M\rightarrow M/G$ into the space
of $G$-orbits on $M$. If the $G$-action is free and proper
in a neighborhood $U$
of $x\in M$, the space $M/G$ is smooth near $p(x)$, and $P_M$ defines
a Poisson structure on $p(U)$.

\section{Generalized moment maps}

In this section we define moment maps for quasi-Poisson actions.
For the actions of a quasi-Poisson Lie group $G_D^{\h}$, 
the space $S=D/G$ is the target of the moment maps. 
We always assume that a moment
map $\mu: M \rightarrow D/G$ is equivariant with respect to the
$G$-action on $M$ and the dressing action on $D/G$. As we show in
Proposition \ref{indep}, the equivariance
condition ensures that the definition is independent of the choice
of an admissible complement $\h \subset \d$.
The following definition of a moment map for a quasi-Poisson action 
of a quasi-Poisson Lie group on a manifold with a bivector 
is a generalization of the notion of equivariant
moment map for the Poisson action of a Poisson Lie group on a Poisson
manifold, defined in \cite{Lu}, 
which itself generalizes the
usual notion of equivariant moment map for a hamiltonian action
of a Lie group on a Poisson manifold.

\subsection{Definition of a moment map}

Let $G_D^{\h}$ be a connected
quasi-Poisson Lie group, and let $x \mapsto x_M$
be the action of $\g$ on $M$ by the infinitesimal generators, defined by
\eqref{inf}, of a
quasi-Poisson action of $G_D^{\h}$ on $(M,P^{\h}_M)$.
We define the map $(P_{M}^{\h})^{\sharp}$ from $1$-forms on $M$ to
vector fields on $M$ by
$<(P_{M}^{\h})^{\sharp}\alpha , \beta > ~ = P_{M}^{\h}(\alpha, \beta)$.

\begin{definition} \label{momentdef}
A map $\mu$ from $M$ to $D/G$, equivariant with respect to the 
action of $G$ on $M$ and to the dressing action of $G$ on $D/G$,
is called a {\it moment map}
for the action of $G_D^{\h}$ on $(M,P^{\h}_M)$ if,
on any open subset $\Omega \subset M$,
\begin{equation} \label{moment}
(P^{\h}_M)^{\sharp}( \mu^* \hat{x}_{\h}) = x_M \ ,                      
\end{equation}
for any $\h$ admissible on $\mu(\Omega)$.
The action of $G_D^{\h}$ on $(M,P^{\h}_M)$ is called 
{\it quasi-hamiltonian } if it admits a moment map. 
A {\it quasi-hamiltonian space} 
is a manifold with a bivector on which a quasi-Poisson Lie
group acts by a quasi-hamiltonian action.
\end{definition}

We observe that, 
in this generalized situation, the target of the moment map is a
quasi-hamiltonian space. 
(The term ``quasi-hamiltonian''
was also used in \cite{AMM1} but in the context of group valued moment 
maps.)
>From \eqref{Smom}, it follows immediately that the 
dressing action of $G$ on $D/G$ is quasi-hamiltonian and has
the identity of $D/G$ as a moment map.
Moreover, let $N$ be any $G$-invariant embedded submanifold of $S=D/G$.
There is a unique bivector $P_{N}$ on $N$ such that the embedding of $N$
into $S$ is a bivector map. Then the quasi-hamiltonian dressing action
of $G$ on $S$ restricts to a quasi-hamiltonian action of $G$ on $N$,
and has the embedding of $N$ into $S = D/G$ as a moment map. In particular,
the orbits of the dressing action of $G$ on $D/G$ which are embedded
submanifolds are quasi-hamiltonian spaces.

\begin{example}
If the Manin quasi-triple $(\d,\g,\h)$ is a Manin triple, then
$\varphi = 0$, so $(G,P_{G})$ is a Poisson Lie group and
$(M,P_{M})$ is a Poisson manifold. In this case, 
the preceding definition of an
equivariant moment map reduces to that given by Lu in \cite{Lu},
that is, the infinitesimal generator of the group action $x_{M}$ is the
image under the Poisson map in $M$ of the pull-back by the moment map of
the right-invariant $1$-form with value $x$ at the identity in the dual 
group $D/G$ of $G$.

In the  particular case of the 
standard quasi-triple, $(T^*G, G, \g^*)$, which corresponds
to the Manin triple $(\g \oplus \g^{*}, \g, \g^{*})$, with $F=0$ and
$\varphi=0$,
the Poisson bivector
$P_{G}$ vanishes and the dual group is the abelian group $\g^{*}$. 
The moment
map then takes values in the vector space
$\g^*$. For any $x \in \g$, the right-invariant 
$1$-form $x^{\rho}$ is the constant form $x$ on $\g^{*}$,
and its pull-back by the moment map
is $d(\mu(x))$, so we recover the usual definition of the moment map
for a hamiltonian
action.
In this case, the orbits of the
dressing action of $G$ on $\g^*$ are the coadjoint orbits familiar from
the usual moment map theory. 
\end{example}

\begin{example} 
\label{diagP}
In the case of $(G\times G, G, \g_-)$,
the quotient space $D/G=G$ is diffeomorphic to the group $G$, and the
dressing action is the action by conjugation. Hence, the 
conjugacy classes in $G$ are quasi-hamiltonian spaces, and the inclusion
is a moment map.
\end{example}

In fact, we do not need to impose the
moment map condition \eqref{moment} for all admissible complements
because conditions \eqref{moment} for different admissible complements
are equivalent.

\begin{proposition} \label{indep}
Let $\h$ and $\h'$ be two complements admissible
at a point $s \in D/G$, and let $m\in M$ be such that $\mu(m)=s$.
Then, at the point $m$, conditions \eqref{moment} for $\h$ and $\h'$
are equivalent, namely 
$ (P^{\h}_M)^{\sharp}( \mu^* \hat{x}_{\h}) 
=  (P^{\h'}_M)^{\sharp}( \mu^* \hat{x}_{\h'})$.
\end{proposition}

\begin{proof}
To prove the proposition, 
we use both \eqref{twisthat} and the above definition of the 
moment map and its equivariance. Thus, at point $m$ such that 
$\mu(m) = s$,
\begin{multline} \nonumber
(P_M^{\h'})^{\sharp}(\mu^* \hat{x}_{\h'})=
(P^{\h}_{M} - t_{M})^{\sharp} (\mu^{*} \hat{x}_{\h'}) = 
(P^{\h}_{M})^{\sharp}(\mu^{*} 
\widehat{\nu x}_{\h}) - t_{M}(\mu^{*}\widehat{\nu x}_{\h})  \\
=(\nu x)_{M} +
((t \circ \tau_{s}) (\nu x))_{M} = ((1_{\g}+ t \circ \tau_{s}) (\nu x))_{M}
= x_{M}.
\end{multline}
Here we have used  the equivariance property of the moment map and
formula \eqref{tau},
which imply that, for any $y$ in $\g$, 
\begin{equation} \nonumber
<\mu^{*}\hat{y}_{\h} , (e_{i})_{M} > (m) 
= ~ <\hat{y}_{\h} , (e_{i})_{S} > (s) = (\tau_{s})_{ik}y^{k},
\end{equation}
when $\mu(m) = s$, and hence
\begin{equation} \label{equiv}
t_{M}(\mu^{*}\hat{y}_{\h}) (m) = - ((t\circ \tau_{s})y)_{M}(m) \ .
\end{equation}
The proposition is therefore proved.
\end{proof}

\subsection{Torus-valued moment maps}

Let us consider the Manin pair 
$(\d, \g)$, where $\d= {\mathfrak u}(1) \oplus {\mathfrak u}(1)$
and $\g = {\mathfrak u}(1)$. Here $\d$ is the direct sum of two copies 
of $\g$, and this is a particular case of Example \ref{gminus}.
There are two group pairs with a compact subgroup corresponding to this
Manin pair, $(S^{1} \times \R, S^{1})$ and $(S^{1}  \times S^{1}, S^{1})$.
The first group pair corresponds to the usual moment map theory.
Let us show that the notion of a moment map for the second group pair 
extends that of an
$S^{1}$-valued moment map in the sense of McDuff \cite{McD} \cite{G}
\cite{W} to the case of a manifold with a possibly degenerate Poisson
bivector. 
A symplectic action of $S^{1}$ on a symplectic manifold $(M, \omega_{M})$ 
is also a quasi-Poisson action of the quasi-Poisson Lie group $S^{1}$ 
defined by the quasi-triple $(S^{1}  \times S^{1}, S^{1} , \g_{-})$,
where $\g_{-}=\Delta_{-}({\mathfrak u} (1))$ (see \ref{gminus}), because
in this case $F=0$ and $\varphi = 0$. 
Here $S = (S^{1}\times S^{1})/S^{1} \simeq S^{1}$.
Since $S^{1}$ is abelian, we see from \eqref{example} that $\hat X =
-  d\alpha$, where $X$ is the generator of $\mathfrak u (1)$ and 
$\alpha$ is the parameter on $S^{1}$. 
A moment map for this quasi-Poisson action is a map, $\mu$,
from $M$ to $S^1$, satisfying 
$  P^{\sharp}_{M}(\mu^{*}d\alpha)=X_M$. When $P_{M}$ is nondegenerate 
with inverse $\omega_{M}$, this condition is equivalent to 
$i_{X_{M}}\omega_{M} = \mu^{*}d\alpha$, where $i$ denotes the interior
product, which is the defining property of an $S^{1}$-valued moment map.

More generally, for $r > 1$, $T^{r}$-valued moment maps (see \cite{G}) 
for symplectic actions of
an $r$-dimensional torus $T^{r} = S^{1} \times \ldots \times S^{1}$
are also moment maps for the quasi-Poisson action of the quasi-Poisson
Lie group $T^{r}$ defined by the quasi-triple $(T^{r} \times T^{r},
T^{r}, \g_{-})$, where $\g$ is the Lie algebra 
${\mathfrak t}^{r} = {\mathfrak u}(1) \oplus \ldots \oplus {\mathfrak u}(1)$,
since in this case also $F = 0$ and $\varphi = 0$, and the preceding
relation is valid for each copy of $S^{1}$.

When $G$ is a compact, connected abelian Lie group,
quasi-Hamiltonian $G$-spaces $(M, \omega_{M})$ in the sense of
\cite{AMM1} are necessarily symplectic. When $M$ is
equipped with the nondegenerate bivector $P_{M}$
defined by $\omega_{M}$, it
is also a quasi-Hamiltonian space for the quasi-Poisson Lie group $G$ 
defined by the quasi-triple $(G \times G, G, \g_{-})$, because in this
case, (\eqref{moment}) coincides with the defining property of the group-valued
moment map of \cite{AMM1}. For nonabelian compact 
Lie groups, it can be shown that
the moment map theory developed here also coincides with 
the moment map theory of \cite {AMM1}.

\subsection{The standard quasi-triple and group-valued moment maps}
\label{standardqt}

We now summarize the case of          
the hamiltonian quasi-Poisson spaces, already in Examples 
\ref{gminus}, \ref{diag}, \ref{ex}, \ref{exP} and \ref{diagP}. 

If $G$ is a connected Lie group with a bi-invariant scalar product,
then there is a well-defined $\Ad$-invariant element $\varphi$   
in $\bigwedge^3\g$, where $\g$ is the Lie algebra of $G$.
If $(K^{\sharp})^{-1}$ 
denotes the isomorphism between $\g^*$ and $\g$ defined by the 
scalar product $ K $ on $\g$, then $\varphi$ satisfies
$$
\varphi(\xi,\eta,\zeta) = {\frac 14}
~<\xi,[(K^{\sharp})^{-1}\eta,(K^{\sharp})^{-1}\zeta]> \ .
$$
The Lie group $G$ diagonally embedded in $G \times G$ defines a group pair
and an associated quasi-triple is 
$(\g \oplus \g, \Delta(\g), {\frac 12}\Delta_{-}(\g))$, with the scalar 
product on $\g \oplus \g$ defined by \eqref{scal}. 
In this quasi-triple, $F=0$ and the element $\varphi$ in  
$\bigwedge^3\g$ is the one defined above.

Let $M$ be a manifold on which the Lie group $G$ acts 
and let $P$ be a bivector field on $M$.
Then $(M,P)$ is a {\it quasi-Poisson space} if $P$ is $G$-invariant and 
\begin{equation}
{\frac 12}[P,P]=\varphi_M \ ,
\end{equation}
where $\varphi_M$ is the field of 3-vectors on $M$ induced 
from $\varphi$ by the infinitesimal action of $\g$ on $M$.

If one identifies $G$ with $S= (G \times G)/G$, then the dressing 
action of $G$ on $S$
is identified with the adjoint action of $G$ on itself.
The $1$-forms $\hat x$ on $G$ defined by the choice of the 
above quasi-triple are such that
\begin{equation} \label{hatMC}
\hat x_g = 2 K((1_{\g}+{\Ad}_g)^{-1}x, \theta_g) \ ,
\end{equation}
for $g \in G$, where $\theta$ is the left Maurer-Cartan form on $G$.

According to Definition 
\ref{momentdef}, $(M,P)$ is called a {\it hamiltonian
quasi-Poisson space} if it is a quasi-Poisson space and if, 
moreover, there exists a moment map for 
the quasi-Poisson action
of $G$ on $M$. 

\begin{proposition}
Let $(M,P)$ be a manifold with a bivector on which the compact simple
Lie group $G$ acts
and which is a quasi-Poisson 
space. Then $(M,P)$ is a quasi-hamiltonian space
if and only if there exists 
a map $\mu: M \rightarrow G$
which is equivariant with respect to the given action 
of $G$ on $M$ and the adjoint action of $G$ on itself,
and which satisfies 
\begin{equation} \label{momentMC}
P^{\sharp}(\mu^*K(x,\theta))={\frac 12} ((1_{\g}+{\Ad}_\mu)x)_M \ ,
\end{equation}
for all $x \in \g$.
\end{proposition}
\begin{proof}
At each point $m \in M$, let us apply the construction 
of Example \ref{ex}, with $g= \mu(m)$, 
to obtain an admissible complement 
by means of which we formulate the definition of the 
moment map. 
The bivector on $M$, 
after the twist $t={\frac {\varepsilon}{2}} \sum_{\alpha \in \Gamma}a_{\alpha}
\wedge b_{\alpha}$, is
$$
P'=P- {\frac {\varepsilon}{2}} \sum_{\alpha \in \Gamma}(a_{\alpha})_{M} 
\wedge (b_{\alpha})_{M} \ .
$$
If $x$ is in the orthogonal complement of the kernel of 
$1_{\g}+{\Ad}_{\mu(m)}$, 
then ${\hat x}'= \hat x$, and $P'(\mu^*{\hat x})=P(\mu^*{\hat x})$. 
Taking into account \eqref{hatMC}, we see that the condition 
\eqref{momentMC} is then 
equivalent to the definition 
$P^{\sharp}_m(\mu^*(\hat x_{\mu(m)}))=x_M(m).$
For $a_{\alpha}, b_{\alpha} \in V_{\alpha}$, 
the moment map conditions are
${P'}^{\sharp}_{M}(\mu^{*}\hat a'_{\alpha})= (a_{\alpha})_{M}, ~~ 
{ P'}^{\sharp}_{M}(\mu^{*}\hat b'_{\alpha})= (b_{\alpha})_{M} \ .$
We replace $\hat a'_{\alpha},
\hat b'_{\alpha}$ by their values found in Example \ref{ex}.
Using the equivariance of the moment map, and the fact that $x^{\rho}= - x^{\lambda}$, 
for $x \in V$, 
we find,
for all $\alpha, \beta \in \Gamma$,
$$
{\frac {\varepsilon}{2}}
\sum_{\alpha \in \Gamma}((a_{\alpha})_{M} 
\wedge (b_{\alpha})_{M})(\mu^*\hat a'_{\beta}) = - (a_{\beta})_M,
$$
$$
{\frac {\varepsilon}{2}} \sum_{\alpha \in \Gamma}((a_{\alpha})_{M} 
\wedge (b_{\alpha})_{M})(\mu^*{\hat b}'_{\beta}) = - (b_{\beta})_M,
$$
Thus the moment map condition, for $x \in V$, reduces to
$$
P^{\sharp}(\mu^{*}K(x, \theta))= 0 \ ,
$$
and therefore it coincides with \eqref{momentMC}.
\end{proof}

\subsection{Generalized foliations of 
quasi-hamiltonian spaces}
We now consider a quasi-hamiltonian $G$-space $(M,P^{\h}_M)$, as in 
Definition \ref{momentdef}, where $G = G^{\h}_D$ is a quasi-Poisson 
Lie group defined by a Manin quasi-triple $(\d,\g,\h)$. 
We wish to show that there is an integrable generalized distribution 
on $M$ defined by the image of any bivector $P^{\h}_M$, 
where $\h$ is admissible.

It follows from the existence of a moment map for the action 
of $G$ on $M$ that, at each point in $M$ where $\h$ is an admissible 
complement, the image of $(P^{\h}_M)^{\sharp}$ 
contains the tangent space to 
the $G$-orbit through this point. 

We further observe that, under a change of admissible complement, when 
$P^{\h}_M$ is modified to $P^{\h'}_M = P^{\h}_M - t_M$ 
(see \eqref{twistP}), the image of $(P^{\h'}_M)^{\sharp}$
coincides with that of $(P^{\h}_M)^{\sharp}$.
In fact, since by the moment map property the image of $t^{\sharp}_M$ is 
contained in the tangent space to the $G$-orbit, then
$\im t^{\sharp}_M \subset \im (P^{\h}_M)^{\sharp}$, and 
therefore $\im (P^{\h'}_M)^{\sharp} = \im (P^{\h}_M)^{\sharp}$.
By symmetry, the two images coincide. We denote by ${\mathcal D}$ the 
generalized distribution thus defined on $M$.

\begin{proposition}
The distribution ${\mathcal D}$ on the quasi-Poisson space $M$
satisfies the Frobenius property, 
$[{\mathcal D}, {\mathcal D}] \subset {\mathcal D}$.
\end{proposition}
\begin{proof}
Let $P = P^{\h}_M$ for an admissible $\h$.
Let $f$ and $g$ be arbitrary functions on $M$, and let 
$[P,P]^{\sharp}(df,dg)$ denote the vector field on $M$ such that
$$
< [P,P]^{\sharp}(df,dg), \alpha > = [P,P](df,dg,\alpha) \ ,
$$
for any $1$-form $\alpha$ on $M$. Then
$$
[P^{\sharp}(df),P^{\sharp}(dg)]- P^{\sharp}d \{ f,g \} 
= - {\frac 12} [P,P]^{\sharp}(df,dg) \ .
$$
Since ${\frac 12} [P,P] = \varphi_M$, and since the tangent space to the $G$-orbit at each point is contained in the distribution 
$\im P^{\sharp}={\mathcal D}$, we conclude that
$[{\mathcal D}, {\mathcal D}] \subset {\mathcal D}$.   
\end{proof}

To conclude that the distribution ${\mathcal D}$ is completely integrable,
it is enough, by the Stefan-Sussmann theorem (see {\it e.g.}, 
\cite{Vaisman}), to show that
the rank of ${\mathcal D}$ is constant along the trajectory of 
any ``hamiltonian vector field'', 
$P^{\sharp}(df), f \in C^{\infty}(M)$.
Let $\mu_t$ be the flow of $P^{\sharp}(df)$ in the neighbourhood
of $m \in M$. (This local one-parameter
group of local diffeomorphims of $M$ depends on $f$.)
For $f, f_{1}, f_{2} \in C^{\infty}(M)$, 
$$
({\mathcal L}_{P^{\sharp}df}P)(df_{1},df_{2})
= {\frac 12} [P,P](df,df_{1},df_{2}) 
= \varphi_M(df,df_{1},df_{2})
=\psi_M(df_{1},df_{2}) \ ,
$$
where $\psi_M$ is a bivector, depending on $f$, with values in the
second exterior power of the tangent space to the $G$-orbits.
It follows that the image of a hamiltonian vector at point $m \in M$ 
under the flow, $\mu_{t}$,
of $P^{\sharp}df$ is an element of the image of $P^{\sharp}$ at point
$\mu_{t}m$. Since the images of a set of linearly independent vectors 
at point $m$ are linearly independent at point $\mu_{t}m$, the dimension
of the image of $P^{\sharp}$ at $\mu_{t}m$  is at least equal to 
that of the
image of $P^{\sharp}$ at $m$. Reversing the argument, we see that the
dimensions of the images at both points are equal.

Thus, there is a well-defined generalized foliation 
(in the sense of Stefan and Sussmann) on $M$, whose leaves 
are nondegenerate quasi-hamiltonian 
spaces, containing the $G$-orbits.
To see that each leaf satisfies the nondegeneracy property, 
we observe that, by the skew-symmetry of $P$, 
the kernel of $P^{\sharp}_m$ 
is the orthogonal of $\im P^{\sharp}_m$.
Therefore, $P^{\sharp}_m$ factorizes through the dual of 
$\im P^{\sharp}_m$ and, by dimension counting, is an isomorphism
onto $\im P^{\sharp}_m$.

\begin{example}
Let the quasi-Poisson space $M$ be $D/G=S$. 
Then the leaves of the generalized foliation 
${\mathcal D}$ are the orbits of the dressing action of $G$ on $S$.

Let, in particular, 
the quasi-Poisson space $M$ be $(G \times G)/G \cong G$ as in 
Section \ref{standardqt}. In this case, 
the orbits of the dressing action are the conjugacy classes of $G$, 
each of which is a nondegenerate quasi-hamiltonian space. 

\end{example}

\subsection{Properties of moment maps}
The generalized moment maps introduced in the previous section possess
properties which resemble the properties of the usual moment maps.

\begin{theorem}
Let $(M,P^{\h}_M,\mu)$ be a quasi-hamiltonian space acted upon by a
connected quasi-Poisson Lie group $G_D^{\h}$. 
Then, on any open set $\Omega\subset M$ such that
$\h$ is admissible on $\mu(\Omega)$, the moment map is a bivector map 
from $(M, P^{\h}_M)$ to $(S,P^{\h}_{S})$,
\begin{equation}
\mu_* P^{\h}_M= P^{\h}_S \ .
\end{equation}

\end{theorem}

\begin{proof}
By the definition of the moment map, for all $x \in \g$, 
\begin{equation} \nonumber
\mu_* x_M= \mu_* (P^{\h}_M)^{\sharp}(\mu^* \hat{x}_{\h}), 
\end{equation}
while, by the characteristic  property of $P_S$,
$x_{S} = (P^{\h}_S)^{\sharp}(\hat{x}_{\h})$.
Thus, we see that 
$\mu_* P^{\h}_M =P^{\h}_{S}$ follows from the equivariance condition,
$\mu_* x_M = x_S$, and from the fact that the $1$-forms $\hat{x}_{\h}$
span the cotangent space to $S$.
\end{proof}

Sometimes it is convenient to require that the bivector
$P_M^{\h}$ be nondegenerate, {\it i.e.}, that the map $(P_M^{\h})^{\sharp}$ 
be an isomorphism.
We next show that for quasi-hamiltonian spaces this condition is 
independent of the particular choice of a complement.

\begin{proposition}
Let $(M,P_M^{\h},\mu)$ be a quasi-hamiltonian space acted upon 
by $G_D^{\h}$.
Let $\h$ and $\h'$
be two isotropic complements of $\g$ admissible on $\mu(\Omega)$, where
$\Omega$ is an open subset of $M$. Then $P_M^{\h}$ is nondegenerate on
$\Omega$ if and only if $P_M^{\h'}$ is nondegenerate on $\Omega$.
\end{proposition}

\begin{proof}
We assume that $P_M^{\h}$ is nondegenerate and we let $\alpha$ be 
a $1$-form in the kernel of $P^{\h'}_M$. Then, 
\begin{equation} \nonumber
(P_M^{\h})^{\sharp}(\alpha)= t_M\alpha.
\end{equation}
Since $\h$ is admissible, there exists $x\in \g$, such
that $t_{M}\alpha = x_{M}$. (If $\alpha = \alpha_{i}(\e^{i})_{M}$, we
can set $x = t^{ik}\alpha_{i}e_{k}$.) By the nondegeneracy of 
$P_M^{\h}$ and the definition of the moment map,
$\alpha= \mu^* \hat{x}_{\h}$. Applying \eqref{equiv}, we obtain 
$x_{M}(m) = - ((t \circ \tau_{s})x)_{M}(m)$, where $s=\mu(m)$.
The equivariance of the moment map implies that
$((1_{\g} +t \circ \tau_{s})x)_{S}=0$, which in turn implies that  
$\tau_{s}(1_{\g} +t \circ \tau_{s}) x = 0$. If both complements
$\h$ and $\h'$ are admissible at $s$, the operator 
$1 + t \circ \tau_{s}$ is invertible, and therefore from 
$(t \circ \tau_{s})(1_{\g} +t \circ \tau_{s}) x = 0$, we obtain
$(t \circ \tau_{s})x = 0$. Since  $x_{M}= - ((t \circ
\tau_{s} )x)_{M}$, the vector field $x_{M}$ vanishes. 
Therefore $\alpha = 0$, and $P^{\h'}_{M}$ is nondegenerate.
\end{proof}

Because of the preceding proposition, it is justified to call
a family of bivectors $P_M^{\h}$ nondegenerate if $P_{M}^{\h}$
is nondegenerate for an admissible $\h$.

\begin{example}
In the case of the standard quasi-triple studied in 
Section \ref{standardqt}, it
follows from \eqref{momentMC} 
that $(M,P,\mu)$ is a nondegenerate quasi-hamiltonian space if and only
if, for each $m \in M$,
$$
{\ker}(P^{\sharp}_m) = \{ \mu^*K(x,\theta) | 
x \in {\ker}(1_{\g}+{\Ad}_{\mu(m)}) \} \ .
$$
\end{example}

Finally, we establish a relation between the Poisson reduction of
Theorem \ref{red} and moment maps. 

\begin{theorem}
Let $(M,P^{\h}_M,\mu)$ be a quasi-hamiltonian space such that the 
bivector
$P_M^{\h}$ is everywhere nondegenerate. 
Assume that $M/G$ is a smooth manifold in a
neighbourhood $U$ of $p(x_{0})$, where $x_{0} \in M$. 
Let $x \in M$ be such that $p(x) \in U$ and $s= \mu(x) \in
D/G$ is a regular value of the
moment map, $\mu$.
Then the symplectic leaf through $p(x)$ in the Poisson manifold $U$ 
is the connected
component of the intersection 
with $U$ of the projection of the manifold $\mu^{-1}(s)$.
\end{theorem}

\begin{proof} The proof is analogous to that in the case of
Poisson actions. See, {\it e.g.}, \cite{Lu} .

Let $x\in M$ be as stated and 
let $y=p(x)$ be its projection in $U \subset M/G$. 
Let $L \subset U \subset M/G$ be the symplectic leaf passing through $y$
in the Poisson manifold $U$.
We choose any complement $\h$ admissible at $s = \mu(x)elm$.
We need to prove that the projection of the tangent space to the
level submanifold, $p_* T_x \mu^{-1}(s)$, coincides with the tangent space
to the symplectic leaf, $T_yL$. 

We denote the Poisson bivector which is locally defined on $M/G$
near $y$ by $Q$.
Let $v$ be a vector in $T_yL$. By definition,
$v=Q^{\sharp}(\alpha)$, where $\alpha$ is a $1$-form in $T^*_y(M/G)$,
and we can also represent $v$ as $v=p_* (P^{\h}_M)^{\sharp}(p^*\alpha)$. 

Vectors $u$ that are 
tangent to the level submanifold $\mu^{-1}(s) \subset M$ are
characterized by the property $\l u, \mu^* \hat{x}_{\h} \r=0$ 
for any $x\in \g$. The 
bivector $P^{\h}_M$ being nondegenerate, any vector $u$ tangent to $M$
at $x$ can be represented as
$u=(P_M^{\h})^{\sharp}(\beta)$ for some $\beta \in T^*_xM$.
Then $u$ is tangent to $\mu^{-1}(s)$ if and only if 
$< \beta , x_M >~=0$, for all $x \in \g$, in other words,
the $1$-form $\beta$ is the inverse image of a $1$-form 
$\alpha \in T^*_y(M/G)$,  $\beta=p^* \alpha$.

We conclude that 
$T_yL= p_* (P^{\h}_M)^{\sharp}(p^* T^*_y(M/G)) = p_* T_x\mu^{-1}(s)$
which proves the theorem.
\end{proof}

Thus, we have extended several of the basic properties of 
hamiltonian group actions on Poisson
manifolds to the case of quasi-hamiltonian actions.

\end{document}